\documentclass[11pt]{article}

\usepackage{latexsym,color,graphicx,amsmath,amsthm,amssymb, subfigure} 

\usepackage{enumitem}

\def\scirc{{\footnotesize \circ}}

\def\R{{\mathbb R}}
\def\S{{\mathbb S}}

\def\C{{\cal C}}

\newtheorem{thm}{Theorem}
\newtheorem{lem}[thm]{Lemma}

\newtheorem{prop}[thm]{Proposition}
\newtheorem*{clm}{Claim}
\begin{document}

\title{On triple intersections of three families of unit circles\thanks{%
Work on this paper by Orit E. Raz and Micha Sharir was supported by Grant 892/13 from the Israel Science Foundation. Work by Micha Sharir was also supported
by Grant 2012/229 from the U.S.--Israel Binational Science Foundation, by the Israeli Centers of Research Excellence (I-CORE) program (Center No.~4/11), and
by the Hermann Minkowski-MINERVA Center for Geometry at Tel Aviv University.  Work by J\'ozsef Solymosi was supported by NSERC, ERC-AdG 321104, and OTKA NK
104183 grants.} \thanks{A preliminary version of this paper appeared in {\em Proc. 30th Annu. ACM Sympos. Comput. Geom.}, 2014, pp. 198--205.}}

\author{
Orit E. Raz\thanks{%
School of Computer Science, Tel Aviv University,
Tel Aviv 69978, Israel.
{\sl oritraz@post.tau.ac.il} }
\and
Micha Sharir\thanks{%
School of Computer Science, Tel Aviv University,
Tel Aviv 69978, Israel.
{\sl michas@post.tau.ac.il} }
\and
J\'ozsef Solymosi\thanks{%
Department of Mathematics,
University of British Columbia,
Vancouver, BC, V6T 1Z4, Canada.
\newline {\sl solymosi@math.ubc.ca}} }
\maketitle

\begin{abstract}
Let $p_1,p_2,p_3$ be three distinct points in the plane, and, for $i=1,2,3$, let $\C_i$ be a family of $n$ unit circles that pass through $p_i$. We address
a conjecture made by Sz\'ekely, and show that the number of points incident to a circle of each family is $O(n^{11/6})$, improving an earlier bound for this
problem due to Elekes, Simonovits, and Szab\'o~\cite{ESSz}. 
The problem is a special instance of a more general problem studied by Elekes and Szab\'o~\cite{ESz} 
(and by Elekes and R\'onyai~\cite{ER00}).
\end{abstract}

\noindent {\bf Keywords.} Combinatorial geometry, incidences, unit circles.

\section{Introduction}

\begin{figure}
\center
\includegraphics[width=0.35\textwidth]{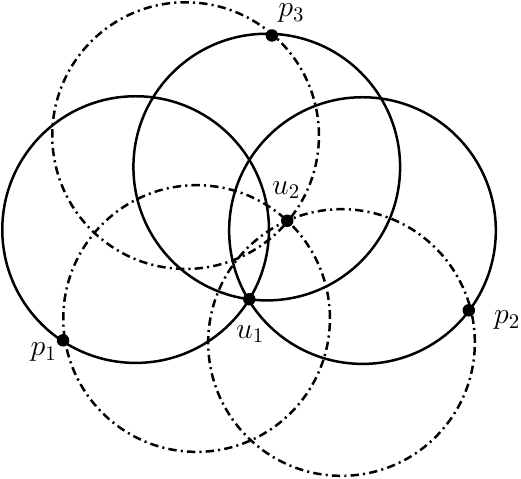}
\caption{Two triple points $u_1$ and $u_2$.}
\label{triples}
\end{figure}

In this paper we re-examine the following problem. Let $p_1,p_2,p_3$ be three 
distinct points in the plane, and, for $i=1,2,3$, let $\C_i$ be a family of
$n$ unit circles that pass through $p_i$. The goal is to obtain an upper bound on 
the number of \emph{triple points}, which are points that are incident to
a circle of each family. See Figures \ref{triples} and \ref{setups}(a) for an illustration. 
Recently, Elekes et al.~\cite{ESSz} have shown that the number
of such points is $O(n^{2-\eta})$, for some constant parameter $\eta>0$ 
(that they did not make concrete); by this they settled a conjecture of Sz\'ekely
(see \cite[Conjecture 3.41]{El}), stipulating that this number should be $o(n^2)$. 
The problem is well motivated in~\cite{ESSz}, because it yields a combinatorial distinction
between unit circles and lines. That is, there exist three families of lines passing 
through three respective points, which determine $\Theta(n^2)$ triple
points, in contrast with the smaller bound of $o(n^{2})$ for unit circles.

Using a different technique, which appears to be simpler than the one in~\cite{ESSz}, 
we establish the following improved bound.
\begin{thm}\label{main}
Let $p_1,p_2,p_3$ be three distinct points in the plane, and, for $i=1,2,3$, let $\C_i$ be a family of $n$ unit circles that pass through $p_i$. Then the
number of points incident to a circle of each family is $O(n^{11/6})$.
\end{thm}

The specific problem studied in this paper can be viewed as a special instance of 
a more general setup, which has been studied by Elekes and
R\'onyai~\cite{ER00} and by Elekes and Szab\'o~\cite{ESz} (see also~\cite{El,ESSz}).
From a high-level point of view the setup is as follows. We have three sets
$A$, $B$, $C$, each of $n$ real numbers, and we have a trivariate real polynomial 
$F$ of some constant degree $d$.  Let $Z(F)$ denote the subset of $A\times
B\times C$ where $F$ vanishes. The claim is that, unless $F$ and $A$, $B$, 
$C$ have some very special structure, $|Z(F)|$ is subquadratic. (For a
simple example where $|Z(F)|$ is quadratic in $n$, consider the case where 
$F(x,y,z)=x+y-z$, and where $A=B=C=\{1,2,\ldots,n\}$.) 

Positive and significant results for this general problem have been obtained by Elekes and R\'onyai~\cite{ER00} and by Elekes and Szab\'o~\cite{ESz}, who
showed that, unless $F$ has a very restricted form, $|Z(F)|$ is indeed subquadratic in $n$. 
For example, in the case where $F$ is of the form $z-f(x,y)$, for some bivariate polynomial $f$,
if $|Z(F)|$ is quadratic in $n$, then $f$ must be of one of the forms $p(q(x)+r(y))$ or $p(q(x)\cdot r(y))$, 
for suitable univariate polynomials $p$, $q$, $r$
(see \cite{ER00} and \cite{El}). 
Related representations, somewhat 
more complicated to state, for a polynomial $F$ with
$|Z(F)|=\Theta(n^2)$, have also been obtained for the general case (see~\cite{ESz} and~\cite{El}).
We have recently studied in \cite{RSS,RSdZ} this specific problem 
and obtained improved bounds for
$|Z(F)|$, when $f$ does not have these special forms; see below.\footnote{The
study in~\cite{RSdZ} has been conducted after the original preparation of this paper;
see a discussion comparing these works in a concluding section.}

The high-level approach used in this paper is similar to those used in several recent works that study problems in combinatorial geometry that are special
instances of this general framework (see Sharir, Sheffer, and Solymosi~\cite{SSS} and Sharir and Solymosi~\cite{SSo}). However, the actual implementations
of this approach in our paper, as well as in the other works just mentioned, are very problem-specific and exploit the special geometric structure of the
relevant problem.

We will later detail the connection of our problem to the setup in \cite{ER00,ESz}. Roughly speaking, for each $\C_i$, its circles have one degree of
freedom, and we parameterize them by a suitable single real parameter. Then, with a proper choice of these parameters, the condition that three circles, one from each family, have a common point can
be expressed by an equation of the form $F(x,y,z)=0$, where $F$ is a real trivariate polynomial, and $x,y,z$ are the parameters representing the three
relevant circles.

In both cases, the specific problem studied in this paper (and the specific problems studied in \cite{SSo,SSS}), and the general one in \cite{ER00,ESz}, the approach is to double count
the number $Q$ of quadruples $(a,p,b,q)$, such that, in our specific context, $a,b$ represent two circles in $\C_1 $, $p,q$ 
represent two circles in $\C_2$, and there exists $z$,
representing a circle through $p_3$, such that $F(a,p,z)=0$ and $F(b,q,z)=0$. (In the general case too, the quadruples 
to be considered are $(a,p,b,q)$ such that $a,b\in A$, $p,q\in B$, and there exists $c\in C$ such that $F(a,p,z)=F(b,q,z)=0$.) A lower bound for $Q$,
in terms of $|Z(f)|$, is easy to obtain via the Cauchy-Schwarz inequality
(see below for details), and an upper bound is obtained by regarding each such quadruple $(a,p,b,q)$ as an {\it
incidence} between the point $(p,q)$, in a suitable parametric plane, and a curve $\gamma_{a,b}$ which is the locus of all points $(p,q)$ that satisfy with
$a,b$ the above conditions. The comparison between the lower and upper bounds yields the asserted upper bound on $|Z(F)|$.

The main issue that arises in bounding the number of incidences is the possibility 
that many curves $\gamma_{a,b}$ overlap each other, in which case the standard
techniques for analyzing point-curve incidences fail. A major part of the analysis in 
this paper is to show that the amount of overlap in our specific problem is bounded. In this
case the standard incidence techniques do apply, and yield a sharp upper bound that 
leads to the aforementioned bound on $|Z(F)|$; see below for details.

In the general problem, the goal is to show that when there is a larger amount of 
overlap between the curves, the polynomial $F$ must have a special form,
as the ones mentioned above and established in \cite{ER00}, and to establish a 
subquadratic upper bound on $|Z(F)|$ when this is not the case. As mentioned,
this indeed has recently been shown in two companion papers, first in \cite{RSS} for the special 
case where $F(x,y,z)=z-f(x,y)$, for any constant-degree bivariate
polynomial $f$, and later in \cite{RSdZ}, for the general trivariate case,
with the same subquadratic bound $O(n^{11/6})$, in both cases, as the one 
established in this paper. In our
problem, though, this part is not needed, and the argument that the overlap is 
bounded is an ad-hoc argument that exploits the geometric
and algebraic structure of the problem.

\section{Unit circles spanned by points on three unit circles} \label{sec:3circs}

\begin{figure}
\subfigure[]{\includegraphics[width=0.35\textwidth]{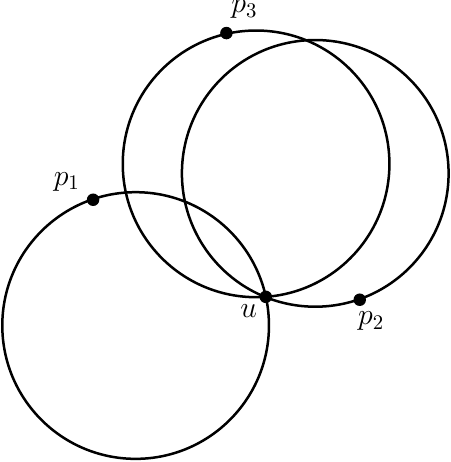}}
\hfill
\subfigure[]{\includegraphics[width=0.35\textwidth]{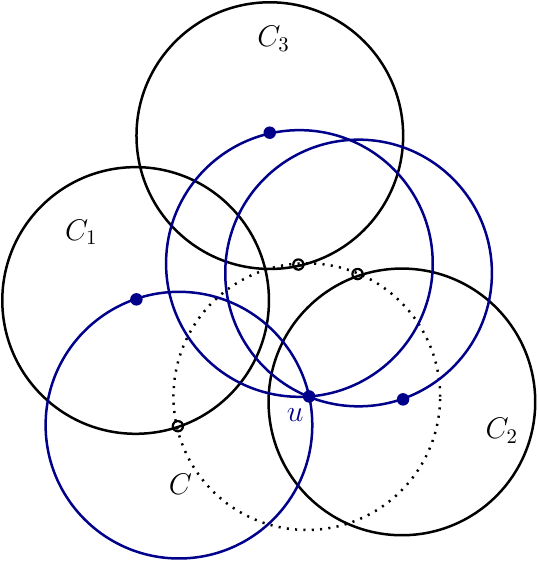}}
\hfill
\caption{(a) A concurrent triple of unit circles and the corresponding triple point $u$.  (b) The triple point $u$ is mapped in the new setup to the unit
circle $C$ centered at $u$. The circles $C_1,C_2,C_3$ are centered at $p_1,p_2,p_3$, respectively, and the hollow points on $C$ are the centers of the three
original circles.}
\label{setups}
\end{figure}

We begin by observing the following equivalent and, in our opinion, more convenient formulation of Theorem~\ref{main}. 
\begin{thm}\label{mainalt}
Let $C_1, C_2, C_3$ be three unit circles in
$\R^2$, and, for each $i=1,2,3$, let $S_i$ be a set of $n$ points lying on $C_i$. Then the number of unit
circles, spanned by triples of points in $S_1\times S_2\times S_3$, is $O(n^{11/6})$.
\end{thm}
The equivalence between this formulation and the one in Theorem~\ref{main} is indeed trivial: 
For each $i$, $S_i$ is the set of centers of the circles of $\C_i$, and the centers of the resulting ``trichromatic'' unit
circles in the new formulation are the triple points in the previous one. 
See Figures~\ref{setups}(a) and \ref{setups}(b) for an illustration of this
connection between the two setups. 
In what follows we prove Theorem~\ref{mainalt}, and stick to the new equivalence formulation.

We note that the condition that three points $p,q,r\in\R^2$ span a unit circle can be expressed as a polynomial equation in their coordinates.
That is, there exists a 6-variate real polynomial $F$ of degree 6, such that $F(p_1,p_2,q_1,q_2,r_1,r_2)=0$ when
$p=(p_1,p_2)$, $q=(q_1,q_2)$, $r=(r_1,r_2)$ span a unit circle. 
Indeed, put
$$
x  = \|p-q\|,\quad X=x^2 ,\quad
y  = \|p-r\|,\quad Y=y^2 , \quad
z  = \|q-r\|,\quad Z=z^2 .
$$
Let $S$ denote the area of the triangle $\Delta pqr$. Then the circumradius $R=1$ of this triangle is given by the formula
$$
1=R=\frac{xyz}{4S} .
$$
The area $S$ can be expressed by Heron's formula, written as
$$
16S^2=(x+y+z)(-x+y+z)(x-y+z)(x+y-z).
$$
That is, we have $(x+y+z)(-x+y+z)(x-y+z)(x+y-z) = x^2y^2z^2$. With some algebraic manipulations, this can be expressed
in terms of the squared distances $X$, $Y$, $Z$, as
\begin{equation} \label{ucirpoly}
X^2 + Y^2 + Z^2 - 2XY - 2XZ - 2YZ + XYZ = 0 .
\end{equation}
The left side of (\ref{ucirpoly}) is the desired polynomial $F$ in the coordinates of $p,q,r$. It is of degree $6$ in these six variables.

For each $i=1,2,3$, each point $p\in S_i$ can be parameterized by (an appropriate algebraic representation of) the orientation $v_p\in\S^1$ of $p$ with
respect to the center $c_i$ of $C_i$ (note that the centers $c_i$ are the points $p_i$ in the original formulation); denote the set of these $n$ orientations as $\Theta_i$. In what follows we will interchangeably use both notations,
referring to a point $p\in S_i$, for $i=1,2,3$, either by its corresponding parameter $v_p\in\Theta_i$, when we want to stress the algebraic nature of the
problem, or as $p$ itself, when geometry is concerned. 

We call a triple $(v_1,v_2,v_3)$, with $v_i\in \Theta_i$, $i=1,2,3$, a \emph{unit triple} if the three 
corresponding points $p_1\in S_1, p_2\in S_2, p_3\in
S_3$ span a unit circle. We use the standard algebraic representation of the $v_i$'s, where we 
replace $v_i$ by $t_i=\tan\frac{v_i}2$, and the corresponding point on $C_i$ then becomes 
$c_i+\left(\frac{1-t_i^2}{1+t_i^2},\frac{2t_i}{1+t_i^2}\right)$. With these representations, the property of being a
unit triple can be expressed by a polynomial equation $f(v_1,v_2,v_3)=0$, obtained by the appropriate 
substitutions into the equation (\ref{ucirpoly}) of $F$. In 
what follows we will refer to the $v_i$'s as orientations, also when 
substituting them in $f$ (the actual substitution should be of the corresponding parameters $\tan\frac{v_i}2$). 
This slightly incorrect treatment is made to simplify the presentation, and has no real consequences in the analysis.

Clearly, $f$ has
constant (and small) degree (which is at most $24$, as is easily checked).
This illustrates how our problem is indeed a special instance 
of the general problem mentioned in the introduction. 

We next argue that, without loss of generality (with a possible re-indexing of the 
input circles and points), we may assume that the points of $S_1$ all lie
in the portion of $C_1$ that lies outside the closed disk circumscribed by $C_2$ 
(this property will become handy for the forthcoming analysis). To see
this, let $D_1,D_2,D_3$ denote the three (closed) unit disks circumscribed by $C_1,C_2,C_3$, 
respectively, and consider the intersection region $K=D_1\cap
D_2\cap D_3$. Assume first that $K$ has a nonempty interior. As is well known, the boundary 
$\partial K$ of $K$ is of the form $c_1\cup
c_2\cup c_3$, where $c_i$ is a single (possibly empty) connected arc of $C_i$, for $i=1,2,3$. 
More generally, in the intersection $K$ of any finite number
of unit disks, each disk contributes at most one connected arc to $\partial K$. Let $C$ be a unit 
circle in the plane, which is not one of $C_1,C_2,C_3$,
and let $D$ denote the disk bounded by $C$. Then, as just mentioned, $C$ contributes a 
single connected arc $c$ to $\partial (K\cap D)$. It follows that $C$
avoids the relative interior of at least one of the arcs $c_1,c_2,c_3$, namely the arc not 
containing any endpoint of $c$. (If one of those arcs is empty,
$C$ trivially misses that arc.) 

It follows that, for every triple $(p_1,p_2,p_3)\in S_1\times S_2\times S_3$ spanning a 
unit circle $C$, at least one of the points $p_1,p_2,p_3$ avoids
$K$, because neither of these points can lie in the interior of $K$, and they lie on $C$, 
which meets $\partial K$ in at most two points. So, for one of the
indices $i_0\in\{1,2,3\}$, and for at least a third of these triples $(p_1,p_2,p_3)$, the 
point $p_{i_0}\in S_{i_0}$ avoids $K$; without loss of generality
assume $i_0=1$. By discarding the other points of $S_1$, we obtain a reduced configuration 
in which the points of $S_1$ lie outside $K$ and the number of
unit triples is at least one third of its original value. That is, each point in (the reduced) $S_1$ 
lies either outside $D_2$ or outside $D_3$. One of
these subsets of $S_1$ participates in at least half the (remaining) unit triples. To recap, by 
removing the points of the other subset, and by re-indexing
if needed, we may assume that all the points of $S_1$ lie outside the disk $D_2$, 
and that the number of unit triples is at least one sixth of the original
number. This reasoning also applies when $K$ is empty or is a singleton 
(with an empty interior), and in fact becomes much simpler in these cases.

We therefore continue the analysis under the assumption that the points of $S_1$ all lie outside $D_2$.

Let $M$ denote the number of unit circles spanned by $S_1\times S_2\times S_3$. 
Our strategy is to double count the quantity $Q$ (mentioned in the
introduction) that we are now going to define. For each $v_3\in\Theta_3$, let $P_{v_3}$ 
denote the set of pairs $(v_1,v_2)\in \Theta_1\times\Theta_2$ such
that $(v_1,v_2,v_3)$ is a unit triple, so $f(v_1,v_2,v_3)=0$. Note that we have 
$M\le \sum_{v_3\in\Theta_3} |P_{v_3}|\le 8M$. Indeed, there are at
most eight triples in $S_1\times S_2\times S_3$ that span the same unit circle 
$C$ ($C$ intersects each of $C_1,C_2,C_3$ in at most two points, and each
triple of points, one from each pair, spans $C$), and clearly, by definition, 
at least one of these triples is counted in $M$.

We now define 
$$
Q := \sum_{v_3\in\Theta_3}|P_{v_3}|^2.
$$ 
The quantity $Q$ may be interpreted as the number of ordered pairs of unit triples of the form
$((v_1,v_2,v_3), (v_1',v_2',v_3))$, with a common third component $v_3$. Using the Cauchy-Schwarz inequality, we have
\begin{equation} \label{eqS}
M \le \sum_{v_3\in\Theta_3} |P_{v_3}| \le
\Big( \sum_{v_3\in\Theta_3}|P_{v_3}|^2 \Big)^{1/2} n^{1/2} = Q^{1/2}n^{1/2} .
\end{equation} 

\paragraph{The curves $\gamma_{a,b}$.}
To obtain an upper bound for $Q$, we use the following approach. Fix two points $a, b\in S_1$, with orientations $v_a,v_b\in \Theta_1$, respectively, and
define $\gamma_{a,b}$ to be the algebraic curve given by the polynomial equation
$$
R(v_x,v_y):={\rm Res}_{v_3}(f(v_a,v_x,v_3), f(v_b,v_y,v_3))=0,
$$
where ${\rm
Res}_{v_3}(f(v_a,v_x,v_3),f(v_b,v_y,v_3))$ is the {\it resultant} of the two polynomials 
$f(v_a,v_x,v_3)$, $f(v_b,v_y,v_3)$ with respect to $v_3$ (which is thus a real bivariate polynomial 
in $v_x,v_y$, independent of $v_3$). By the properties of the resultant, the curve $\gamma_{a,b}$ 
contains all points $(v_x,v_y)$, with corresponding points $x,y\in C_2$, for which there exists $v_3$ (not necessarily in $\Theta_3$)
such that
\begin{align} \label{eqgam}
f(v_a,v_x,v_3)  &= 0,\\
f(v_b,v_y,v_3)  &= 0 \nonumber ;
\end{align}
for more details see, e.g., Cox et al.~\cite{CLO}.
However, $\gamma_{a,b}$ might contain points $(v_x,v_y)$ where there is no real point $z\in
C_3$ that spans unit circles with both pairs $(a,x)$, $(b,y)$. 
In general, the curve $R(v_x,v_y)=0$ is partitioned into a constant number of connected arcs of two kinds: {\em real} arcs, 
over which (\ref{eqgam}) has a real solution $v_3$, and {\em non-real} arcs, over which 
there are no such real solutions. 
We refer to the endpoints of these arcs as {\em transition points}. 
We will analyze and handle these points later on.\footnote{In both real and non-real arcs, 
we only consider real values for $v_\xi, v_\eta$. Only $z$ can assume non-real values for 
points on non-real arcs. In other words, ignoring its geometric interpretation, $\gamma_{a,b}$ is a real curve.}

Let $\Pi$ denote the set $\Theta_2\times\Theta_2$, represented as a set of points in the 
above parametric plane, let $\Gamma'$ denote the (multi-)set of the
curves $\gamma_{a,b}$, and let $I'=I'(\Pi,\Gamma')$ denote the number of incidences 
between the curves of $\Gamma'$ and the points of $\Pi$.  

Note that, for any fixed $v_3\in\Theta_3$ and for any ordered pair of pairs $(v_a,v_c)$, 
$(v_b,v_d)$ in $P_{v_3}$, we have $(v_c,v_d)\in \gamma_{a,b}$ and
$(v_d,v_c)\in \gamma_{b,a}$. It follows that the number $I'$ of point-curve  
incidences is at least $\frac14 \sum_{v_3\in \Theta_3}|P_{v_3}|^2$. Indeed,
there can be at most four values of $v_3$ that give rise to the same incidence 
(any of the pairs $(v_a,v_c),(v_b,v_d)$, say $(v_a,v_c)$, defines at most two
unit circles that pass through the two corresponding points, and each of these 
circles can intersect $C_3$ in at most two points), and only those values
among them that belong to $\Theta_3$ are reflected in the above sum; also, the 
fact that each pair of pairs in $P_{v_3}$ generate two incidences is
``neutralized'' by the fact that the same two incidences are generated for each of 
the two orderings of the pairs.  That is, we have $Q \le 4I'$, so it
suffices to obtain an upper bound for $I'$.  

The number of incidences between curves of $\Gamma'$ that are of the form 
$\gamma_{a,a}$, with $a\in S_1$, and the points of $\Pi$, is $O(n^2)$. Indeed, let
$a\in S_1$ and consider the curve $\gamma_{a,a}$. For $c\in S_2$, there exist at 
most two unit circles that pass through $a$ and $c$, and these circles form
at most four intersection points $z$ with $C_3$. Then for each such intersection 
point $z$, there exist at most two unit circles that pass through $a$
and $z$, which form at most four intersection points $d$ with $C_2$. Thus there 
are at most 16 values $d$ for which $(v_c,v_d)\in \gamma_{a,a,}$. It
follows that, for each $a\in S_1$, the curve $\gamma_{a,a}$ is incident to $O(n)$ 
points of $\Pi$, and hence the total number of incidences that curves of
this form contribute is $O(n^2)$. 

Therefore, letting $\Gamma\subset\Gamma'$ be the (multi-)set of the curves 
$\gamma_{a,b}$, with $a\neq b\in S_1$, and letting $I=I(\Pi,\Gamma)$ denote the
number of incidences between the curves of $\Gamma$ and the points of $\Pi$, 
we get $I'\le I+O(n^2)$, and thus $Q \le 4I+O(n^2)$. 

We reiterate that $I$ (and $I'$) might include many irrelevant incidences, first, 
because the corresponding parameter $v_3$ does not belong to $\Theta_3$,
and second, because $v_3$ is not real (the incidence occurs on a non-real arc of 
the relevant curve $\gamma_{a,b}$). Still, an upper bound on the
(overestimate) $I$ suffices for our purpose.

Hence the problem is reduced to obtaining an upper bound on $I$. This is an instance 
of a fairly standard point-curve incidence problem, which can be
tackled using the well established machinery, such as the incidence bound of 
Pach and Sharir~\cite{PS98}, or, more fundamentally, the crossing-lemma
technique of Sz\'ekely~\cite{Sz97} (on which the analysis in \cite{PS98} is based). 
However, to apply this machinery, it is essential that the curves of
$\Gamma$ have a constant bound on their multiplicity. More precisely, we need to 
know that no more than $O(1)$ curves of $\Gamma$ can share a common
irreducible component. In more detail, while the points of $\Pi$ are clearly distinct, 
there might be potentially many pairs of curves of $\Gamma$ that coincide or
overlap in a common irreducible component, in which case the aforementioned incidence-bounding 
techniques break down. Fortunately, this can be controlled through the following
key proposition. (Recall that this arises as a key issue when applying this approach to 
the general setup of Elekes and R\'onyai~\cite{ER00}, as manifested
in the companion paper \cite{RSS}, and the more general setup of Elekes and Szab\'o~\cite{ESz},
as in the more recent study~\cite{RSdZ}.)

\begin{prop} \label{pr:ov}
There exists a subset $U\subset S_1\times S_1$ of size at most $O(n)$ such that the following holds.
For any irreducible component $\gamma'$, 
there are at most $O(1)$ pairs $(a,b)\in (S_1\times S_1)\setminus U$ such that $\gamma'$ contains a portion of a real arc of $\gamma_{a,b}$.
\end{prop}

The proof of the proposition is given in Section~\ref{se:ov}. This allows us to derive an upper bound on the number of incidences, 
given in the following proposition. 
\begin{prop}\label{prop:incid}
Let $\Gamma$ and $\Pi$ be as above. 
Then the number $I$ of incidences between $\Gamma$ and $\Pi$ is $O(|\Gamma|^{2/3}|\Pi|^{2/3}+|\Gamma|+|\Pi|)$.
\end{prop}

The proof of the proposition is given in Section~\ref{se:incid}. Since $|\Pi|, |\Gamma| = O(n^2)$, it follows that in this case
$I=O\left(n^{8/3}\right)$ and thus, recalling that $Q\le 4I+O(n^2)$, $Q=O(n^{8/3})$ too, so we get $M=O\left(n^{11/6}\right)$.
This completes the proof of Theorem~\ref{mainalt}. $\hfill\square$

\subsection{\bf Properties of the curves $\gamma_{a,b}$}
In this section we provide a detailed analysis of the structure and properties of the 
curves $\gamma_{a,b}$, from both algebraic and geometric perspectives.

\paragraph{Explicit construction of the third point of a unit triple.}
We slightly change the notation temporarily, and let $a=(a_1,a_2)$ be a point on $C_1$, and $x=(x_1,x_2)$ be a point on $C_2$.
We derive below an explicit expression for a point $z=(z_1,z_2)$ on $C_3$ 
such that $(a,x,z)$ span a unit circle. This procedure will be used repeatedly in the forthcoming analysis.
We note that the procedure consists of two similar substeps. Later on, we will formally break it into these
substeps, and use them as the primitive building blocks for the analysis.

\begin{figure}
\subfigure[]{\includegraphics[width=0.42\textwidth]{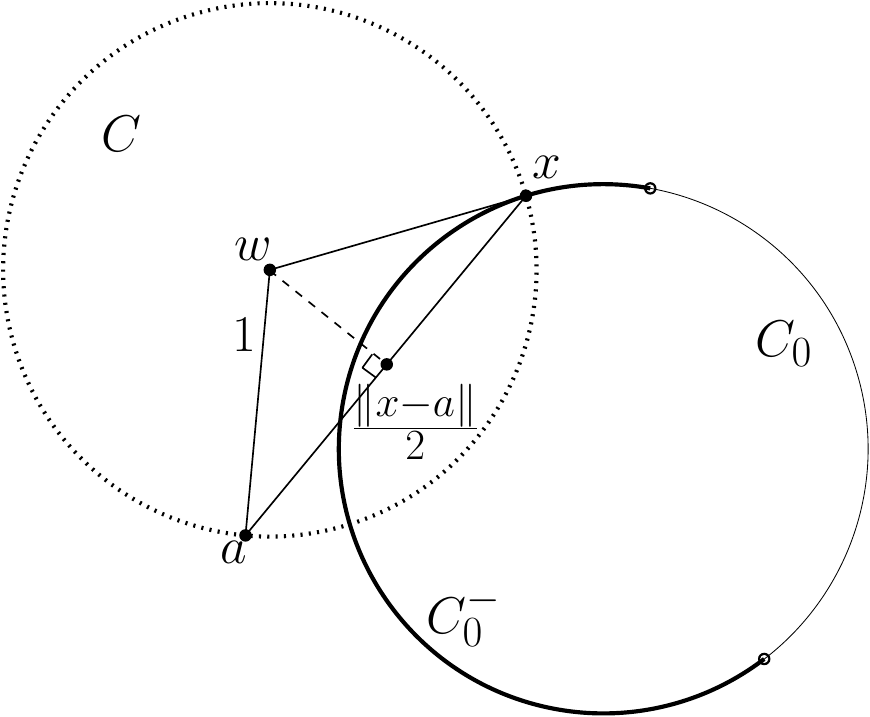}}
\hfill
\subfigure[]{\includegraphics[width=0.45\textwidth]{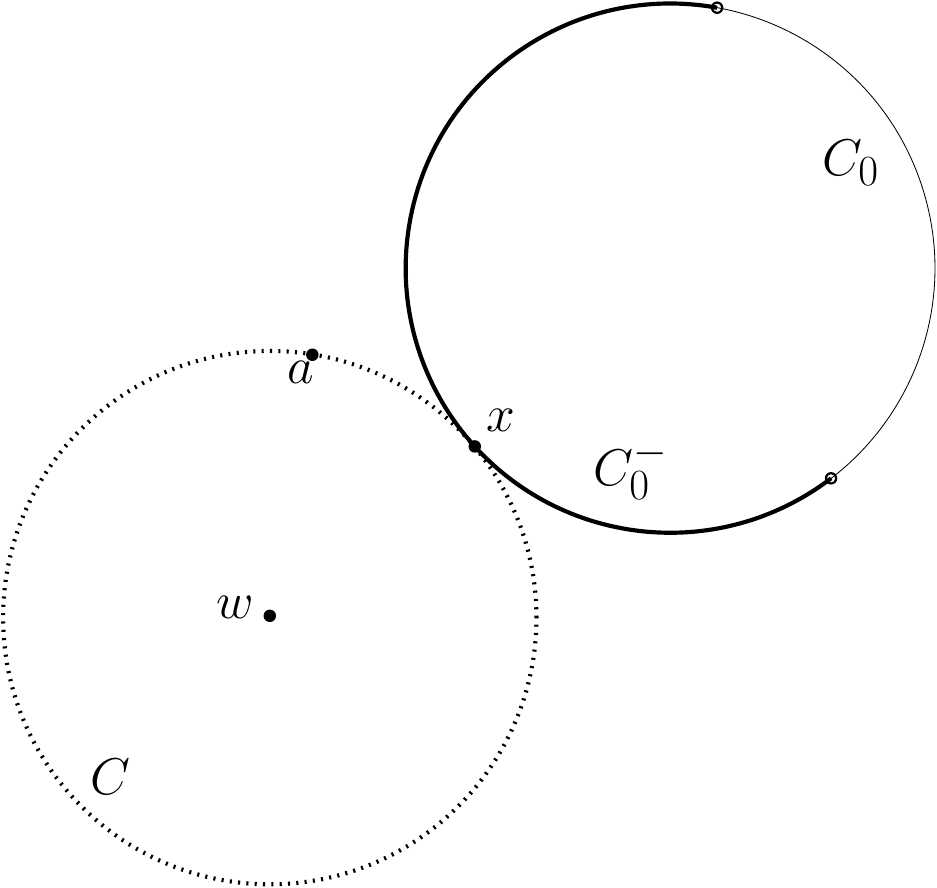}}
\caption{(a) Computing the center $w$ of a unit circle $C$ passing through $a$ and $x$. Here $x$ varies along $C_0$ and $w$ exists when $x$ lies in $C_0^-$ (the highlighted arc). (b) The exceptional situation in Lemma~\ref{derivs}, where $C$ is tangent to $C_0$ at $x$.}
\label{explicit}
\end{figure}

In full generality, let $a$ and $x$ be any pair of points in the plane, 
where we think of $a$ as fixed and of $x$ as a variable.
Let $C$ be a unit circle that passes through $a$ and $x$. The center $w$ of $C$ is the point
\begin{equation}\label{center}
w=(w_1,w_2)=\left(\frac{a_1+x_1}{2},\frac{a_2+x_2}{2}\right)\pm s\left(-\frac{x_2-a_2}{2},\frac{x_1-a_1}{2}\right),
\end{equation}
where $\displaystyle s=\frac{\sqrt{1-\frac{\|x-a\|^2}{4}}}{\frac{\|x-a\|}{2}}=\sqrt{\frac{4}{\|x-a\|^2}-1}$; 
see Figure~\ref{explicit}(a). There can be zero, one, or two real solutions for $w$.
Denote this doubly-valued function as $w=\varphi_a(x)$.
In this notation, both $x$ and $w$ have two degrees of freedom.
However, in our application $x$ will be assumed to lie on some unit circle, $C_0$
so it will have only one degree of freedom, and then
$w$ also has one degree of freedom.
We capture this extra constraint by using the 
modified notation $w=\varphi_{a,C_0}(x)$,
where now $x$ is constrained to lie on $C_0$.

Let $c_3=(q_1,q_2)$ be the center of $C_3$. Then, similar to (\ref{center}), 
the point $z$ is given by $z=\varphi_{c_3,C_a}(w)$. That is,
\begin{equation}\label{z}
z=(z_1,z_2) = \left(\frac{w_1+q_1}{2},\frac{w_2+q_2}{2}\right)\pm r\left(-\frac{q_2-w_2}{2},\frac{q_1-w_1}{2}\right),
\end{equation}
where $\displaystyle r=\sqrt{\frac{4}{\|w-c_3\|^2}-1}$; again, for a given $w$, there 
can be zero, one, or two real solutions for $z$. Thus, the number of real values
of the combined expression for $z$, obtained by substituting (\ref{center}) into 
(\ref{z}), is between $0$ and $4$. 

Symmetric constructions give explicit expressions for the first or second point of a unit triple
in terms of the two other points.

\noindent{\bf Remark.} When $a$ and $x$ are given, the equation $F(a,x,z)=0$, as given in (\ref{ucirpoly}), is 
of degree 4 in the coordinates of $z$. The combination of (\ref{center}) and (\ref{z}) is in fact an explicit solution of 
this equation by radicals (with the standard, one-dimensional representations of the points involved).

\paragraph{The primitive steps of the procedure.}
A couple of additional remarks are in order.
First, similar to what has just been remarked, the roles of $a$ and $x$ in defining $w$ are 
essentially identical. The above notation is supposed to signify that $a$ is considered as a 
fixed parameter and $x$ as a variable (along its circle $C_0$).
Second, $\varphi$ is in general 2-valued, unless $\|a-x\|=2$ (it has complex values 
when $\|a-x\|>2$). In what follows, we will always trace a single (real) branch of 
such a $\varphi$ (over which $\|a-x\|<2$), but stop when $\|a-x\|$ becomes $2$.

\begin{figure}
\subfigure[]{\includegraphics[width=0.4\textwidth]{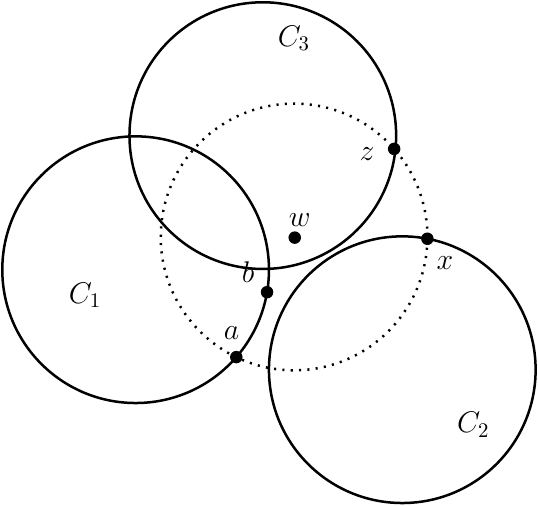}}
\hfill
\subfigure[]{\includegraphics[width=0.4\textwidth]{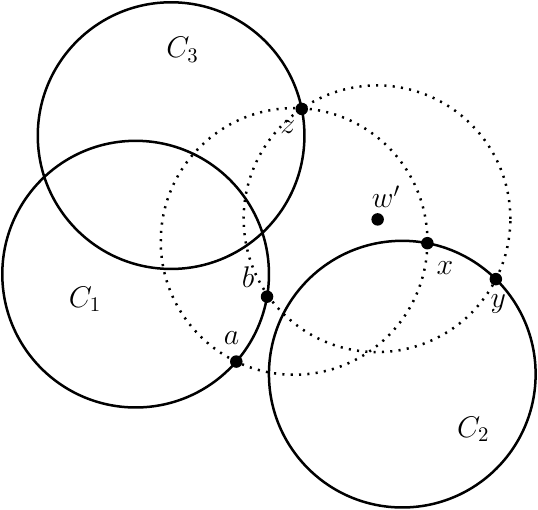}}
\caption{(a) Steps (i) and (ii) of the construction. (b) Steps (iii) and (iv) of the construction. }
\label{const}
\end{figure}

\paragraph{Explicit construction of points on a curve $\gamma_{a,b}$.}
Let $\gamma_{a,b}$ be one of our curves. The preceding construction immediately leads to the following
4-step procedure for constructing points on $\gamma_{a,b}$. Specifically, given $a,b\in C_1$ 
and a point $x\in C_2$, we compute point(s) $y\in C_2$ such that $(v_x,v_y)\in\gamma_{a,b}$, as follows. 
\begin{description}
\item{(i)} We start with $a$ and $x$, and construct the point(s) $w=\varphi_{a,C_2}(x)$, each of which is
the center of a unit circle passing through $a$ and $x$. See Figure~\ref{const}(a). We fix one of these points $w$. 
(We terminate the procedure, with failure, when there are no real solutions; 
this also applies to each of the following steps.)
\item{(ii)} We then construct the point(s) $z=\varphi_{c_3,C_a}(w)\in C_3$, where $C_a$ is
the unit circle centered at $a$. By construction, we have $F(a,x,z)=0$; 
again, there are (at most) two choices for $z$ and we fix one of them. See Figure~\ref{const}(a).
\item{(iii)} We now construct the center(s) $w'=\varphi_{b,C_3}(z)$ of the unit circle(s) passing 
through $b$ and $z$. See Figure~\ref{const}(b). Fix $w'$ to be one of these centers.
\item{(iv)} Finally, we obtain the desired point(s) $y=\varphi_{c_2, C_b}(w')\in C_2$, 
where $C_b$ is the unit circle centered at $b$. See Figure~\ref{const}(b). 
\end{description}
(Note that the circles appearing in the four applications of the $\varphi$-functions, 
namely, $C_2$, $C_a$, $C_3$, and $C_b$, are indeed fixed, regarding $a$ and $b$ themselves as fixed.)

The following lemma shows that the functions $\varphi$ are well-behaved, in the sense
made precise below, unless certain degenerate situations arise.
In the lemma, the derivative of $w=\varphi_{a,C_0}(x)$ should have the expected
interpretation. Concretely, let $v_x$ (resp., $v_w$) denote the orientation of $x$ 
(resp., of $w$) with respect to the center of $C_0$ (resp., $a$). Interpret $w=\varphi_{a,C_0}(x)$ 
as the corresponding functional relationship between $v_x$ and $v_w$, which, with a slight
abuse of notation, we also write as $v_w = \varphi_{a,C_0}(v_x)$. Then $\varphi'_{a,C_0}(x)$ is
simply the corresponding derivative of $\varphi_{a,C_0}(v_x)$ at the respective orientation $v_x$.

\begin{lem} \label{derivs}
Let $a$ be a fixed point, $C_0$ a fixed (unit) circle, and assume $a\notin C_0$. Put $C_0^-:= \{x\in C_0 \mid \|a-x\| < 2\}$. 
Fix a parameterization of $C_0^-$, $x(t)=(x_1(t),x_2(t))$, for $t\in I\subseteq\S^1$, such that each of the 
functions $x_1(t), x_2(t)$ is an analytic function on $I$. 
Then each branch of the function $\varphi_{a,C_0}(x)$ is analytic, and 
\begin{equation}\label{deriv}
\varphi'_{a,C_0}(x) = \frac{(w-x)\cdot\tau_x}{(w-x)\cdot\tau_w} ,
\end{equation}
where $w=\varphi_{a,C_0}(x)$, and $\tau_x$ and $\tau_w$ denote the unit 
tangent vectors to $C_0$ at $x$ and to $C_a$ at $w$, respectively.
In particular, $\varphi_{a,C_0}(x)$ has non-zero derivative, 
at each point $x\in C_0^-$ for which the unit circle centered at $\varphi_{a,C_0}(x)$
is not tangent to $C_0$ at $x$.
\end{lem}
\noindent{\bf Proof.}
See Figure~\ref{explicit}(a) for the general layout, and Figure~\ref{explicit}(b) for the exceptional situation. 
First note that our assumption, combined with the explicit expression (\ref{explicit}), implies that $\varphi_{a,C_0}$ is analytic on $I$. 
Fix a point $x\in C_0^-$ and put $w=\varphi_{a,C_0}(x)$. 
Let $v_x$ and $v_w$ denote the corresponding orientations, let $x'$ be the point with $v_{x'}=v_x + \Delta v_x$,
for a small increment $\Delta v_x$, and put $w'=\varphi_{a,C_0}(x')$ and $\Delta v_w = v_{w'}-v_w$. 
Clearly, since $\varphi_{a,C_0}(x)$ is continuous over $I$, $\Delta v_w$ is also small when $\Delta v_x$ is small. 

Put $\Delta x = x'-x$ and $\Delta w = w'-w$ (note that these are vector displacements, whereas
$\Delta v_x$, $\Delta v_w$ are scalars, and we have $\|\Delta x\|<|\Delta v_x|$, $\|\Delta w\|<|\Delta v_w|$). We have
\begin{align*}
1 & = \|w'-x'\|^2 = \|w+\Delta w - x-\Delta x\|^2 \\
& = \|w-x\|^2 + 2(w-x)\cdot(\Delta w-\Delta x) + o(|\Delta v_x|+|\Delta v_w|) \\
& = 1 + 2(w-x)\cdot(\Delta w-\Delta x) + o(|\Delta v_x|+|\Delta v_w|) .
\end{align*}
Let $\tau_x$ and $\tau_w$ denote the unit tangent vectors to $C_0$ at $x$ and to $C_a$ at $w$,
respectively. We then have
\begin{align*}
\|\Delta x - (\Delta v_x) \tau_x \| & = o(|\Delta v_x|) \\
\|\Delta w - (\Delta v_w) \tau_w \| & = o(|\Delta v_w|) .
\end{align*}
We thus have
$$
(w-x)\cdot((\Delta v_w) \tau_w - (\Delta v_x) \tau_x) = o(|\Delta v_x|+|\Delta v_w|) .
$$
That is,
$$
((w-x)\cdot\tau_w) \Delta v_w = ((w-x)\cdot\tau_x) \Delta v_x + o(|\Delta v_x|+|\Delta v_w|) .
$$
The vector $w-x$ cannot be orthogonal to $\tau_w$. Indeed, this is possible only when either $a$ 
coincides with $x$, which has been ruled out in the lemma, or when $\|a-x\|=2$, which is also excluded. 
Hence the coefficient of $\Delta v_w$ is nonzero,
so we have
$$
\frac{\Delta v_w}{ \Delta v_x} = \frac{(w-x)\cdot\tau_x}{(w-x)\cdot\tau_w}  
+ o\left(1+ \frac{|\Delta v_w|}{|\Delta v_x|}\right) ,
$$
so in the limit we get
$$
\varphi'_{a,C_0}(x) = \frac{(w-x)\cdot\tau_x}{(w-x)\cdot\tau_w} ,
$$
which is well defined and nonzero unless $w-x$ is orthogonal to $\tau_x$ 
(that is, the unit circle centered at $w$ is tangent to $C_0$), 
as asserted in the lemma.
$\Box$

Note that each of the numerator and the denominator of (\ref{deriv}) can become $0$ as $x$ varies 
along $C_0$: The numerator becomes $0$ when $w$ is at distance $2$ from the center of $C_0$ (the exceptional case depicted in Figure~\ref{explicit}(b)), and the
denominator becomes $0$ when $x$ is at distance $2$ from $a$ (the endpoints of $C_0^-$, if they exist). 
Note that both situations are ruled out in the lemma.
These kinds of degeneracy will play a central role in what follows.

Let $a,b$ be a pair of points on $C_1$ (outside the disk $D_2$), and let $(v_x,v_y)$ be a point on $\gamma_{a,b}$, for a respective pair of points $x,y\in C_2$.
Consider the 4-step procedure that produces $y$ from $x$, as described above, and write the outputs of its steps as 
$w=\varphi_{a,C_2}(x)$, $z=\varphi_{c_3,C_a}(w)$, $w'=\varphi_{b,C_3}(z)$, and $y=\varphi_{c_2,C_b}(w')$. In each step
we pick an appropriate branch of the relevant function, and assume that the degeneracies ruled out
in Lemma~\ref{derivs} (and summarized in the preceding paragraph) do not arise (in suitable neighborhoods of the four respective points) in any of these 
four steps. With these notation, assumptions, and conventions, we obtain the composite function
$$
v_y = \Phi_{a,b}(v_x) := \varphi_{c_2,C_b} \scirc \varphi_{b,C_3} \scirc \varphi_{c_3,C_a} \scirc \varphi_{a,C_2} (v_x) ,
$$
which is well defined and analytic in a suitable neighborhood $N$ of $x$, its graph
over $N$ is a portion of $\gamma_{a,b}$, and $(v_x,v_y)$ is not a local $x$-extremum or $y$-extremum 
of that portion. The non-extremality properties are consequences of the chain rule, combined with 
a four-way application of (\ref{deriv}) in the proof of Lemma~\ref{derivs}. (See below for a more explicit repetition of this
argument.)

\subsection{\bf Proof of Proposition~\ref{pr:ov}}\label{se:ov}

Let $\gamma'$ be an irreducible component of potentially many curves $\gamma_{a,b}$. The strategy is to 
identify (few) points along $\gamma'$, 
from which we can reconstruct all the values of $a$ and $b$ in only a constant number of ways. 

Fix a generic point $q_0=(v_{x_0},v_{y_0})$ on $\gamma'$ that is non-singular and non-extremal for $\gamma'$, 
and which does not lie on any other irreducible component of any other curve. Let $\Xi(q_0)$
denote the subset of all pairs $(a,b)$ so that $\gamma_{a,b}$ contains $\gamma'$, and $q_0$ 
lies on a real arc of $\gamma_{a,b}\cap \gamma'$. The preceding discussion implies that, for each such curve
$\gamma_{a,b}$, there is a unique way to define the corresponding function $\Phi_{a,b}$ in a 
suitable neighborhood of $q_0$ (which may depend on $a$ and $b$), and the graph of each of 
these functions near $q_0$ (in the intersection of all these neighborhoods) is $\gamma'$ itself.

Now trace $\gamma'$, along with all these functions, from $q_0$ in, say, increasing $v_x$-direction,
and stop at the first point $(v_\xi,v_\eta)$ at which one of the assumptions made in 
Lemma~\ref{derivs} is violated, for some pair $(a,b)\in \Xi(q_0)$, for one of the corresponding functions
$w=\varphi_{a,C_2}(\xi)$, $z=\varphi_{c_3,C_a}(w)$, $w'=\varphi_{b,C_3}(z)$, or $\eta=\varphi_{c_2,C_b}(w')$. 
The definition of the curves $\gamma_{a,b}$, and the property that all the points of $S_1$ lie outside the disk $D_2$, 
 imply that $\gamma'$ does indeed contain such a point $(v_\xi,v_\eta)$ (because the open disk of radius 2 centered at $a$ intersects 
but does not contain the circle $C_2$).
As we will shortly argue, there are no transition points of any curve $\gamma_{a,b}$, for $(a,b)\in\Xi(q_0)$, 
between $q_0$ and $(v_\xi,v_\eta)$, so $(v_\xi,v_\eta)$ still belongs to the same real subarcs of all 
the curves with pairs in $\Xi(q_0)$ (possibly being a transition point of some of these arcs).

Applying explicitly the chain rule, we have, for an arbitrary point $(v_x,v_y)$ on the relatively 
open portion $(\gamma')^-$ of $\gamma'$ between $q_0$ and $(v_\xi,v_\eta)$, with corresponding points $x,y\in C_2$,
$$
\Phi'_{a,b}(x) = \varphi'_{c_2,C_b}(w') \varphi'_{b,C_3}(z) \varphi'_{c_3,C_a}(w) \varphi'_{a,C_2} (x) ,
$$
where $w$, $z$, and $w'$ are as defined above. By (\ref{deriv}), this is a product of four fractions, and
the above assumption about $(v_\xi,v_\eta)$ means that the numerator or denominator of at least 
one of these fractions is zero at $\xi$, but they all remain nonzero before reaching $\xi$. 
Call $(a,b)$ an \emph{ultra-degenerate} pair if (at least) one numerator and one denominator of $\Phi'_{a,b}$ vanish
simultaneously at $\xi$. 

We will shortly show that ultra-degenerate pairs are sparse,
and handle them all via a global argument that does not depend on $\gamma'$. For the time
being, we ignore all such pairs. (The exceptional set $U$ in the proposition will be the set of
these ultra-degenerate pairs.)

In other words, excluding ultra-degenerate pairs, $\Phi'_{a,b}(v_x)$, which is equal to the slope of the tangent at the corresponding
point on $\gamma'$ (assuming that point to be non-singular), becomes zero or tends to $\infty$ at $v_\xi$, so $\gamma'$ has a (one-sided)
horizontal or vertical tangent at $(v_\xi,v_\eta)$. 

For technical reasons that will become clearer later on, we trace $\gamma'$ from $q_0$ in both increasing and decreasing $v_x$-direction.

It is important to note that the curve $v_y=\Phi_{a,b}(v_x)$ does indeed trace $\gamma'$ between
$q_0=(v_{x_0},\Phi_{a,b}(v_{x_0}))$ and $(v_\xi,v_\eta)$. Indeed, let $h(t,s)=0$ be the (irreducible) 
polynomial equation defining $\gamma'$. Let $t_1\in(v_{x_0},v_\xi)$ be such that the graph of $\Phi_{a,b}$, 
restricted to the subinterval $[v_{x_0},t_1]$, is an arc of $\gamma'$, but this does not hold for $[v_{x_0},t_1^+]$, 
for any $t_1^+>t_1$, arbitrarily close to $t_1$. We thus have $h(t, \Phi_{a,b}(t))=0$, 
for every $t\in[v_{x_0},t_1]$. By Lemma~\ref{derivs}, the function 
$\Phi_{a,b}(t)$ is analytic in some sufficiently small neighborhood $N$ of $t_1$ inside $(v_{x_0},v_\xi)$. 
Thus, letting  $H(t)=h(t,\Phi_{a,b}(t))$, and since $h$ is a polynomial, we have that $H$ is analytic and $H(t)=0$, 
for every $t\in [v_{x_0},t_1]\cap N$. In particular, the derivatives of $H$ at $t=t_1$, of any order, are all zero 
(because $H$ is identically zero in a one-sided neighborhood of $t_1$). Thus, since $H$ is analytic, the 
Taylor series of $H$ at $t_1$ is identically zero, which means that $H$ is identically zero in some suitable 
neighborhood $N'\subset N$ of $t_1$ (this time on both sides of $t_1$). 
In other words, the graph of $\Phi_{a,b}$ continues to coincide with $\gamma'$ on the other side 
of $(t_1,\Phi_{a,b}(t_1))$ too, contradicting our assumption on $t_1$. 

We note that the preceding argument also holds when the tracing of $\gamma'$ encounters a singular point $(v_x,v_y)$ of $\gamma'$ (before reaching $(v_\xi,v_\eta)$). Even if two branches of $\gamma'$ meet tangentially at $(v_x,v_y)$, the graph of $\Phi_{a,b}$ remains well defined, and follows a unique branch of $\gamma'$, on either side of $(v_x,v_y)$.  

\paragraph{Transition points.}
Recall that a point $(v_\xi,v_\eta)$ is a transition point of the 
curve $\gamma_{a,b}$ if it connects a real arc  and a non-real arc of $\gamma_{a,b}$.
We argue that at a transition point we must have one of the degeneracies ruled out in 
Lemma~\ref{derivs}, for one of the four $\varphi$-functions.
Specifically, let $(v_\xi,v_\eta)$ be a point on $\gamma'$ which is a transition point along a containing 
curve $\gamma_{a,b}$. Each of the four functions $\varphi_{a,C_2}$, $\varphi_{c_3,C_a}$, 
$\varphi_{b,C_3}$, $\varphi_{c_2,C_b}$, whose composition yields $\Phi_{a,b}$,
involves a square root (with a fixed sign). As long as none of these roots vanishes, the 
functions continue to be defined (as real functions) and produce real values, so the two
unit circles that are spanned by the corresponding triples $(a,\xi,z)$, $(b,\eta,z)$,
for a suitable point $z\in C_3$, are such that $z$ is real and the circles are real too,
and this continues to hold in a suitable neighborhood of $(\xi,\eta)$. Since this does
not occur at a transition point, one of the square roots has to vanish at $(\xi,\eta)$, 
and when this occurs one of the degenerate conditions in the lemma occurs for the corresponding function
(where the denominator of one of the fractions in (\ref{deriv}) vanishes).
This establishes the promised claim.

\paragraph{Recap.} The preceding analysis leads to the following overall treatment of $\gamma'$. We partition $\gamma'$ into maximal connected subarcs, each delimited by points with horizontal or vertical tangency (and does not contain any such point in its relative interior). 
Since $\gamma'$ has constant degree, there are only $O(1)$ subarcs of this kind. For each of the containing curves $\gamma_{a,b}$, each such subarc $\gamma''$ is fully contained either in a real arc of $\gamma_{a,b}$ or in a non-real arc of $\gamma_{a,b}$, and at least one subarc $\gamma''$ is contained in a real arc of $\gamma_{a,b}$. 

In the next step of the analysis, we show that, for each of the $O(1)$ locally extremal points $(v_\xi,v_\eta)\in\gamma'$, there are only $O(1)$ pairs $(a,b)$, for which $(v_\xi,v_\eta)$ is a (possibly delimiting) point of a real arc of $\gamma_{a,b}$. Altogether, we conclude that $\gamma'$ can be an irreducible component of only $O(1)$ curves $\gamma_{a,b}$ (with the additional requirements that $\gamma'$ overlaps at least one real arc of $\gamma_{a,b}$ and that the pair $(a,b)$ is not ultra-degenerate).


\paragraph{The possible geometric scenarios near a locally extremal point of $\gamma'$.}
Let $(v_\xi,v_\eta)$ be a locally $x$-extremal point of $\gamma'$ (locally $y$-extremal points will be handled in a fully symmetric manner; see below), and let $\xi,\eta$ 
be the points in $C_2$ with orientations $v_\xi,v_\eta$, respectively. Let $a,b\in C_1$ 
be a fixed pair for which at least one of the subarcs of $\gamma'$ delimited by 
$(v_\xi,v_\eta)$ is a (portion of a) real arc of $\gamma_{a,b}$, and assume that $(a,b)$ is not an ultra-degenerate pair. Here we do not know $a$ and $b$, 
and our goal is to reconstruct them from $v_\xi, v_\eta$. This is done as follows.

We first note that, as argued earlier, the $x$-extremality of $(v_\xi,v_\eta)$ means that $\Phi_{a,b}'(v_\xi)=\infty$ for every such pair $(a,b)$. That is, one of the denominators in the four expressions, as in (\ref{deriv}), vanishes. We thus have the following four respective situations.

\begin{figure}
\subfigure[]{\includegraphics[width=0.4\textwidth]{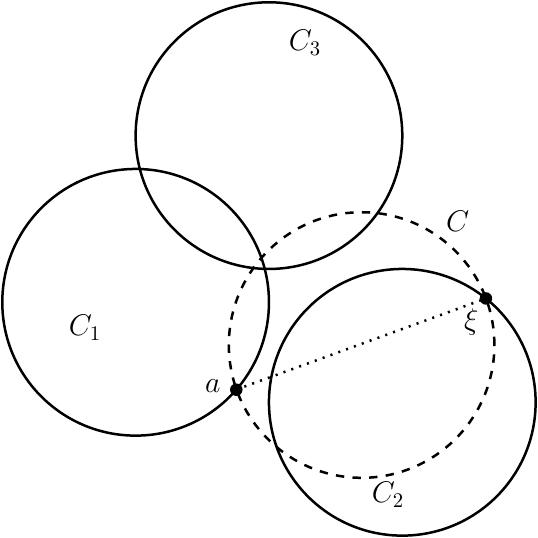}}
\hfill
\subfigure[]{\includegraphics[width=0.4\textwidth]{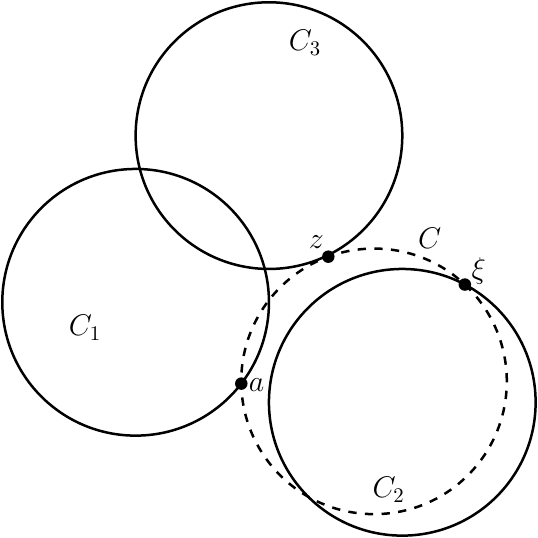}}
\caption{(a) The situation in Case (i) of the reconstruction ($a\xi$ is a diameter of $C$). 
(b) The situation in Case (ii) ($C$ is tangent to $C_3$ at $z$).}
\label{fail12}
\end{figure}

\paragraph{Case (i) $\|a-\xi\|=2$.}
In this case we can reconstruct $a$ in at most two possible
ways, as an intersection point of $C_1$ with the circle of radius 2 centered at $\xi$; see Figure~\ref{fail12}(a). We can then retrieve the corresponding point $z\in C_3$ in two possible ways, as an intersection
point of $C_3$ with the (unique) unit circle that passes through $a$ and $\xi$. Since $\eta$ is also given, we can compute $b$, as one of the intersection
points of $C_1$ with one of the at most two unit circles that pass through $z$ and $\eta$. Altogether, there are (at most) two ways to choose $a$, two for
$z$, and four for $b$, so in the present case we can reconstruct $(a,b)$ in at most 16 possible ways.

\paragraph{Case (ii) $\|c_3-w\|=2$.}
In this case, there is a unit circle that passes through $a$ and $\xi$ and is tangent to $C_3$ at $z$; see Figure~\ref{fail12}(b). Hence $z$ is a tangency point of
$C_3$ with one of the at most two unit circles that are incident to $\xi$ and tangent to $C_3$. This allows us to reconstruct $a$, as an intersection point
of $C_1$ with one of these two unit circles. We then retrieve $b$ from $z$ and $\eta$ as in the preceding case. Altogether, there are (at most) two ways to choose $z$, two for
$a$, and four for $b$, so here too we can reconstruct $(a,b)$ in at most 16 possible ways.

\begin{figure}
\subfigure[]{\includegraphics[width=0.4\textwidth]{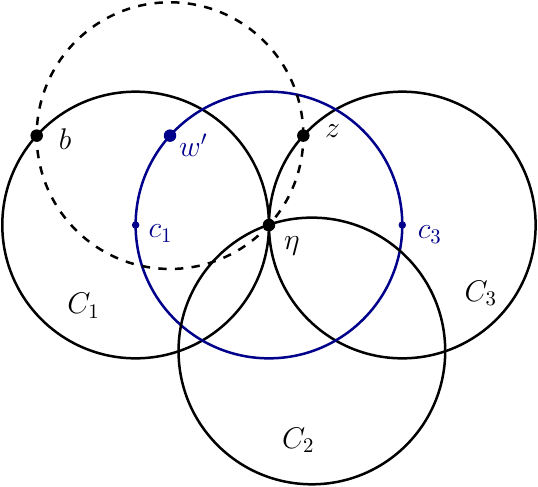}}
\hfill
\subfigure[]{\includegraphics[width=0.4\textwidth]{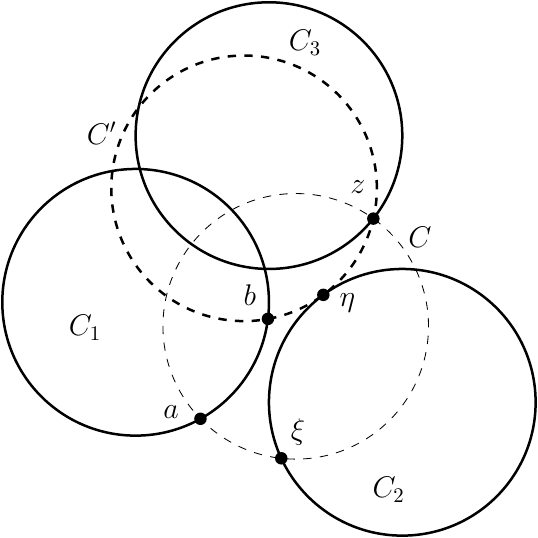}}
\caption{(a) The exceptional situation in the claim in Case (iii) ($C_1$ and $C_3$ are tangent at $\eta$). 
(b) The situation in Case (iv) ($C'$ is tangent to $C_2$ at $\eta$).}
\label{fail34}
\end{figure}

\paragraph{Case (iii) $\|b-z\|=2$.}
This is geometrically more challenging to analyze, because the points $z,b$, at distance 2 apart, are both unknown. See Figure~\ref{fail34}(a). We handle this case by observing that the lengths of
the edges of the quadrilateral $R=c_1bzc_3$ are fixed---they are 1,2,1 and $|c_1c_3|$, respectively, but this does not determine $R$, because it can
flex (with one degree of freedom) about its fixed edge $c_1c_3$. As $R$ flexes, the midpoint $w'$ of $bz$ traces an algebraic curve $\tau$ of some constant
degree $d$. Note that the unit circle that passes through $b,\eta$ and $z$, has its center $w'$ (which is the midpoint of $bz$) on $\tau$. Since the point $\eta$
is known, we can find $w'$, by computing the intersection points of $\tau$ with the unit circle $C_\eta$ centered at $\eta$, and then retrieve $b$, as the
intersection point of $C_1$ with the unit circle centered at $w'$. We claim that, in general, there are at most $2d$ intersection points of $\tau$ with
$C_\eta$, and hence at most a constant number of ways to reconstruct $b$. Indeed, if this were not the case, then, by B\'ezout's theorem (see, e.g.,
\cite{CLO}), $\tau$ would have to contain $C_\eta$ as one of its components. This situation is controlled by the following simple claim.

\begin{clm}
The curve $\tau$ does not contain any unit circle as one of its components, unless $C_1$ and $C_3$ are tangent to each other, in which case $\tau$ does indeed contain
the unit circle centered at the point of tangency.
\end{clm}
\noindent{\bf Proof.}
See Figure \ref{fail34}(a)  for the exceptional situation in the claim. Let $C$ be a unit circle, centered at a point $c$, 
such that $C\subseteq \tau$. By the construction of $\tau$, every point $p\in C$ is the midpoint of a
segment whose endpoints lie on $C_1$ and $C_3$, respectively. This implies, in particular, that $C$ is contained in $K:={\rm conv}(C_1\cup C_3)$, the convex
hull of $C_1\cup C_3$, and since the three circles $C, C_1,C_3$ are of the same radius, it follows that the center $c$ of $C$ lies on $c_1c_3$. Moreover, as is easily checked (cf.~Figure \ref{fail34}(a)), in this case $c$ must be the midpoint of $c_1c_3$, and we must have $|c_1c_3|=2$, and then $C_1$ and $C_3$ are tangent to each
other at $c$, implying that $C_1$, $C_3$, and $C$ must indeed be of the exceptional kind stated in the claim. $\Box$

To complete the analysis, we recall that the only problematic case is when the unit circle centered at $\eta$ is contained in $\tau$. By the claim, $\eta$
must be the midpoint of $c_1c_3$, so in particular $C_2$ is not tangent to $C_1$ (or else it would coincide with $C_3$). Note that in this special situation, the points of $S_1$, with the
possible exception of $\eta$, clearly lie outside the disk circumscribed by $C_3$. We can discard the tangency point $\eta$ of $C_1$ and $C_3$ from $S_1$,
if needed, losing at most $O(n)$ unit triples spanned by the original sets $S_1,S_2,S_3$. We now restart the whole analysis, switching the roles of $C_2$
and $C_3$, and are now guaranteed that the exceptional situation described in the above claim does not occur.

We can then retrieve $z$, as an intersection point of $C_3$ with a unit circle that passes through $b$ and $\eta$, and, since $\xi$ is also given, compute
$a$, as one of the intersection points of $C_1$ with one of the unit circles that pass through $z$ and $\xi$. 

\paragraph{Case (iv) $\|c_2-w'\|=2$.}
In this case, depicted in Figure~\ref{fail34}(b), the unit circle that passes through $b,\eta$ and $z$ is tangent to $C_2$ at $\eta$. Hence $b$ is one of
the intersection points of $C_1$ and the unique unit circle (externally) tangent to $C_2$ at $\eta$, and then $a$ can be reconstructed from $b,\xi$ and $\eta$, as in the previous cases.

\paragraph{Handling $y$-extremal points of $\gamma'$.}
The cases where $(v_\xi,v_\eta)$ is a $y$-extremal point of $\gamma'$ are handled in a fully symmetric manner, as follows. Let $\hat\gamma'$ denote the ``transpose" of $\gamma'$, that is, 
$$
\hat\gamma'=\{(v_q,v_p)\mid (v_p,v_q)\in\gamma'\}.
$$
Clearly, $\gamma'$ is an irreducible component of $\gamma_{a,b}$ if and only if $\hat\gamma'$ is an irreducible component of $\gamma_{b,a}$. Moreover, $(v_\xi,v_\eta)$ is a locally $y$-extremal point of $\gamma'$ if and only if $(v_\eta,v_\xi)$ is a locally $x$-extremal point of $\hat\gamma'$. We can therefore apply the preceding analysis, essentially verbatim, to $\hat\gamma'$ and the ``transposed" containing curves $\gamma_{b,a}$, and conclude that, from each locally $y$-extremal point $(v_\xi,v_\eta)$ of $\gamma'$, there are only $O(1)$ curves $\gamma_{a,b}$ such that $\gamma'\subset\gamma_{a,b}$ and $(v_\xi,v_\eta)$ lies on a real arc of $\gamma_{a,b}$.

\paragraph{Ultra-degenerate pairs.}
Let $U$ denote the set of all ultra-degenerate pairs; recall that a pair $(a,b)$ is said to be {\em ultra-degenerate} if there exists a point $(\xi,\eta)$
at which at least one numerator and one denominator of the four fractions that define
$\Phi'_{a,b}$ vanish simultaneously at $\xi$. The reason for singling out these pairs
is that it is not clear what happens to the slope of the tangent to $\gamma_{a,b}$ (that is, to $\gamma'$) at this point.

Fortunately, the overall size of $U$, over all possible components $\gamma'$, is only $O(n)$. 
The somewhat tedious case analysis that establishes this claim is given in Appendix~\ref{app:ultra}. 
Since each such curve has only $O(n)$ incidences with the points
of $S_2\times S_2$ (an easy property which is a special case of the Schwartz--Zippel Lemma~\cite{Sch80,Zi89}), we get a total of $O(n^2)$ incidences that
correspond to such pairs. This is a small bound, subsumed by the overall bound on $Q$
that we derive.

Combining all the steps of the analysis, and excluding the $O(n)$ pairs in $U$, we finally conclude that $\gamma'$ overlap (a portion of a real arc of) $\gamma_{a,b}$, for at most a constant number of pairs $(a,b)$.
This completes the proof of Proposition~\ref{pr:ov}. $\Box$


\subsection{\bf Proof of Proposition \ref{prop:incid}}\label{se:incid}
We apply Sz\'ekely's technique~\cite{Sz97}, which is based on the crossing lemma (see also Pach and Agarwal~\cite{PA95}). As noted, this is also the
approach used in \cite{PS98}, but the possible overlap of curves requires some extra (and more explicit) care in the application of the technique. A similar
argument is given in the companion paper~\cite{RSS}, but we repeat it here to make the paper more self-contained.

Throughout this subsection, we completely ignore curves $\gamma_{a,b}$ for which $(a,b)$ is an ultra-degenerate pair; as argued in Section~\ref{se:ov}, these curves contribute only $O(n^2)$ to the incidence bound.

We begin by constructing a plane embedding of a multigraph $G$, whose vertices are the points of $\Pi$, and each of whose edges connects a pair $\pi_1=(\xi_1,\eta_1)$,
$\pi_2=(\xi_2,\eta_2)$ of points that lie on the same curve $\gamma_{a,b}$ and are consecutive along (some connected component of) $\gamma_{a,b}$; the edge
is drawn along the portion of the curve between the points. One edge for each such curve (connecting $\pi_1$ and $\pi_2$) is generated, even when the curves
coincide or overlap. Thus there might potentially be many edges of $G$ connecting the same pair of points, whose drawings coincide. Nevertheless, by
Proposition~\ref{pr:ov}, this number is at most $O(1)$.

In spite of this control on the number of mutually overlapping (or, rather, coinciding) edges, we still face the potential problem that the edge
multiplicity in $G$ (over all curves, overlapping or not, that connect the same pair of vertices) may not be bounded (by a constant). More concretely, we
want to avoid edges $(\pi_1,\pi_2)$ whose multiplicity exceeds $d^2$, where here $d\le 576$ 
denotes the degree of the curves of $\Gamma$. 

To follow this strategy, we pass to a dual parametric plane, in which the roles of $\Theta_1$ and $\Theta_2$ are interchanged, so curves $\gamma_{a,b}$ of
$\Gamma$ become dual points $(v_a,v_b)$, and points $(v_\xi,v_\eta)$ of $\Pi$ become dual curves $\gamma_{\xi,\eta}^*$, defined as the locus of all points
$(v_x,v_y)$, each corresponding to a pair of points $x,y\in C_1$, for which there exists $v_3$ (not necessarily in $\Theta_3$) such that (compare with (\ref{eqgam}))
\begin{align}
f(v_x,v_\xi,v_3) & = 0 \\
f(v_y,v_\eta,v_3) & = 0 ; \nonumber 
\end{align}
we denote by $\Gamma^*$ and $\Pi^*$ the sets of the dual points and of the dual curves, respectively. Clearly, we have $(v_\xi,v_\eta)\in \gamma_{a,b}$ if and only if $(v_a,v_b)\in\gamma_{\xi,\eta}^*$. We recall our assumption that, for
$(v_a,v_b)\in\Gamma^*$, the corresponding points $a$ and $b$ lie on the portion of $C_1$ which is outside $D_2$. In view of this, we can ignore irreducible
components of curves of $\Pi^*$ which contain only ``irrelevant" points $(v_a,v_b)$, that is, only points $(v_a,v_b)$ for which one of $a$ or $b$ lies in $D_2$. 

\begin{clm}
Let $\pi_0=(\xi_0,\eta_0)\in\Pi$ be a vertex of $G$. Then there exist at most 
$d^3-1$ vertices $\pi\in\Pi$, such that the number of curves of $\Gamma$ 
that pass through both $\pi_0,\pi$ is larger than $d^2$.
\end{clm}

\noindent{\bf Proof.}
For contradiction, assume there exist $\pi_i=(\xi_i,\eta_i)$, $i=1,\ldots,d^3$, 
such that, for each $i$, there are at least $d^2+1$ curves of $\Gamma$ passing through 
both $\pi_0,\pi_i$.
By construction (and using duality), every curve of $\Gamma$ connecting $\pi_0$ and $\pi_i$ 
corresponds to a dual point of $\Gamma^*$ that lies on both dual curves
$\gamma_{\xi_0,\eta_0}^*$ and $\gamma_{\xi_i,\eta_i}^*$. Thus, by the assumption 
on the pair $(\pi_0,\pi_i)$, for each $i=1,\ldots, d^3$, the
curves $\gamma_{\xi_0,\eta_0}^*$ and $\gamma_{\xi_i,\eta_i}^*$ have at least $d^2+1$ 
points in common, and hence, by B\'ezout's theorem (see, e.g.,
\cite{CLO}), the two curves share a common irreducible component. Note that 
$\gamma_{\xi_0,\eta_0}^*$, having degree $d$, has at most $d$ irreducible
components, and thus, by the pigeonhole principle, there exists an irreducible 
component $\gamma_0^*$ of $\gamma_{\xi_0,\eta_0}^*$ that is shared by at
least $d^2$ curves $\gamma_{\xi_i,\eta_i}^*$; by reindexing, if needed, assume 
these are $\gamma_{\xi_i,\eta_i}^*$, $i=1,\ldots,d^2$. 

Let $\kappa=O(1)$ be the multiplicity bound obtained in Proposition~\ref{pr:ov}. Let $(v_{a_j},v_{b_j})$, $j=0,\ldots, \kappa d$, be $\kappa d+1$ (distinct) points on $\gamma_0^*$, having the property that, for each $j=0,\ldots, \kappa d$, the points
$a_j,b_j$ lie on the portion of $C_1$ which is outside $D_2$ (but not necessarily in $S_1$); as already noted, we may assume that $\gamma_0^*$ contains at
least one such point, but then, by continuity, it contains infinitely many such points. We have $(v_{a_j},v_{b_j})\in \gamma_{\xi_i,\eta_i}^*$, or, by
duality, $(v_{\xi_i},v_{\eta_i})\in \gamma_{a_j,b_j}$, for every $i=1,\ldots,d^2,\; j=0,\ldots, \kappa d$. Similarly, we also have
$(v_{a_j},v_{b_j})\in\gamma^*_{\xi_0,\eta_0}$ and so $(v_{\xi_0},v_{\eta_0})\in\gamma_{a_j,b_j}$, for $j=0,\ldots,\kappa d$.  Using B\'ezout's theorem once
again, we have that $\gamma_{a_0,b_0}$ and $\gamma_{a_j,b_j}$, which intersect in at least $d^2+1$ points, share a common irreducible component, for each
$j=1,\ldots,\kappa d$. Since $\gamma_{a_0,b_0}$ is of degree $d$, and thus has at most $d$ irreducible components, we conclude that there exists an irreducible
component of $\gamma_{a_0,b_0}$ that is shared by $\kappa$ other curves $\gamma_{a_j,b_j}$. This however contradicts 
Proposition \ref{pr:ov}, and hence the
claim follows. $\Box$

Consider a point $\pi_1$ and one of its bad neighbors\footnote{We make the pessimistic assumption 
that they are (consecutive) neighbors along all
these curves, which of course does not have to be the case in general.} $\pi_2$. Let $\gamma_{a,b}$ 
be one of the curves along which $\pi_1$ and $\pi_2$ are
neighbors. Then, rather than connecting $\pi_1$ to $\pi_2$ along $\gamma_{a,b}$, we continue along 
the curve from $\pi_1$ past $\pi_2$ until we reach a good
point for $\pi_1$ (i.e., a point such that the number of curves of $\Gamma$ 
that pass through both $\pi_1,\pi$ is at most $d^2$),
and then connect $\pi_1$ to that point (along $\gamma_{a,b}$). We skip over at most $d^3-1$ points in the process, but now, having
applied this ``stretching'' to each pair of bad neighbors, each of the modified edges has multiplicity at most $2d^2$ (the factor 2 comes from the fact that
a new edge $e$ can be obtained by stretching an original edge from either endpoint of $e$).

Note that this edge stretching does not always succeed: It will fail when the connected component $\gamma'$ of $\gamma_{a,b}$ along which we connect the
points contains fewer than $d^3+1$ points of $\Pi$, or when there are fewer than $d^3-1$ points of $\Pi$ between $\pi_1,\pi_2$,
and the ``end" of $\gamma'$ (recall the constraint in the definition of the curves $\gamma_{a,b}$). Still, the number of new edges in $G$ is at least $I(\Pi,\Gamma)-\lambda|\Gamma|$, for a suitable constant $\lambda$, where the
term $\lambda|\Gamma|$ accounts for missing edges on connected components of the curves, for the reasons just discussed. By what have just been argued, the
number of edges lost on any single component is at most $O(d^3)$.

The final ingredient needed for this technique is an upper bound on the number of 
crossings between the (new) edges of $G=(V,E)$.  
Before doing so, we note that the way in which the edges are drawn, some of them may pass 
through vertices of $G$ (other than their endpoints), which was not allowed in Sz\'ekely's original work~\cite{Sz97}.
We resolve this issue by slightly perturbing the edges, so that they are drawn slightly off the curves.
Each crossing between a pair of new edges, before this perturbation, is essentially
a crossing between two curves of $\Gamma$.  Even though the two curves might overlap in a common 
irreducible component (where they have infinitely many
intersection points, none of which is a crossing), the number of proper 
crossings between them is $O(d^2)=O(1)$, as follows, for example, from the
Milnor--Thom theorem (see \cite{Mil64,Th65}), or from B\'ezout's theorem. 
Finally, because of the way the drawn edges have been stretched, the edges, even
those drawn along the same original curve $\gamma_{a,b}$, may now overlap one 
another, or, in the actual pertubed way in which they are drawn, may ``tangentially
cross'' one another. In this case a crossing between two curves may be claimed by more than
one pair of (stretched) crossing edges, and the tangential crossings need to be accounted for too.  
Nevertheless, since no edge straddles more 
than $d^3-1$ points, the number of pairs that claim a specific crossing
is still a constant (that depends on $d$), and so is the number of their tangential crossings.  
Hence, we conclude that the total number 
of edge crossings in $G$ is $O(|\Gamma|^2)$.

We can now continue by applying the crossing lemma argument, exactly as done by Sz\'ekely and in other works (e.g., see \cite{PS98,Sz97}).
The crossing lemma asserts that $\frac{|E|^3}{|V|^2}\le c{\rm Cr}(G)$,
for a suitable constant $c$ (that now depends on $d$, to account for the possible overlap between edges),
provided that $|E|\ge c'|V|$, for another constant $c'$ (that also depends on $d$).
Combining the two possibilities, of a large $|E|$ and a small $|E|$, and using the fact that 
$|E|\ge I(\Pi, \Gamma)-\lambda|\Gamma|$, $|V|=|\Pi|$, and $|{\rm Cr}(G)=O(|\Gamma|^2)$, we obtain
$$
I(\Pi,\Gamma) = O\left( |\Pi|^{2/3}|\Gamma|^{2/3} + |\Pi| + |\Gamma| \right) ,
$$
with the constant of proportionality depending on $d$. This completes the proof
of Proposition~\ref{prop:incid}.
$\Box$


\section{Conclusion}
We do not know whether the bound in Theorems~\ref{main} and \ref{mainalt} is 
tight in the worst case, and suspect that it is not. Resolving this question is a major problem
for further research, especially since it arises in all the related specific and general problems 
in \cite{ER00, ESSz, ESz, RSS, RSdZ, SSo, SSS}.

This problem is clearly only one special instance of several 
related problems, in which we have three sets $S_1,S_2,S_3$ of points,
each contained in some curve, and we want to bound the number of triples in 
$S_1\times S_2\times S_3$ that satisfy some property (that can be specified by
polynomial equation, such as spanning a unit circle). As a simple example, 
consider the case where each $S_i$ 
is contained in some respective line $\ell_i$, for $i=1,2,3$, and the property
is that the triple span a triangle of unit area. In a companion paper \cite{RS15} 
we show that in this case the number of triples can be $\Theta(n^2)$, but
the bound is likely to drop when the sets $S_i$ are contained in other curves.

It is constructive to compare the study in this paper with the
general setup studied in Elekes and Szab\'o~\cite{ESz} and in the recent paper \cite{RSdZ}
(prepared after the original submission of this paper).
In principle, the results in this paper could be interpreted as a special case
of the analysis in \cite{ESz,RSdZ}. That is, we have a trivariate polynomial $F$ over
a Cartesian product $S_1\times S_2\times S_3$ of three sets of $n$ real numbers each,
and the number of the unit circles (or triple points) under consideration is the number 
of zeros of $F$ in $S_1\times S_2\times S_3$.
The results in \cite{RSdZ} assert that the number of such zeros is $O(n^{11/6})$,
{\it unless} $F$ has a special form.
Although a general procedure for such a test can be provided (see, e.g., a discussion in \cite{RS15}),
its concrete execution is extremely complicated --- at the moment we do not see any 
reasonable way to carry it out.
We note that this is a general issue, that one will face in any specific application 
of the general theory of Elekes and Szab\'o.

\paragraph{Acknowledgment.} We would like to thank an anonymous referee, 
whose helpful and constructive comments helped us a lot to improve the paper. 
Part of this research was performed while the authors were visiting the Institute for Pure 
and Applied Mathematics (IPAM), which is supported by the National Science Foundation. 
The authors deeply appreciate the stimulating environment and facilities provided by IPAM.

\appendix

\section{Ultra-degenerate pairs}\label{app:ultra}

Recall that a pair $(a,b)$ is said to be {\em ultra-degenerate} if there exists a point $(\xi,\eta)$
at which at least one numerator and one denominator of the four fractions that define
$\Phi'_{a,b}$ vanish simultaneously at $\xi$. The reason for singling out these pairs
is that it is not clear what happens to the slope of the tangent to $\gamma_{a,b}$ at this point.

\begin{figure}
\subfigure[]{\includegraphics[width=0.49\textwidth]{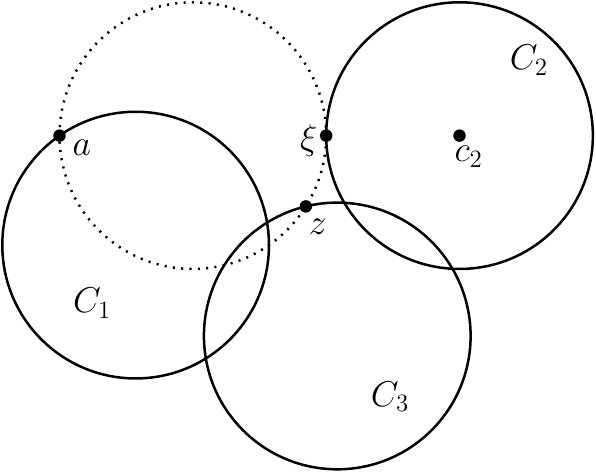}}
\hfill
\subfigure[]{\includegraphics[width=0.45\textwidth]{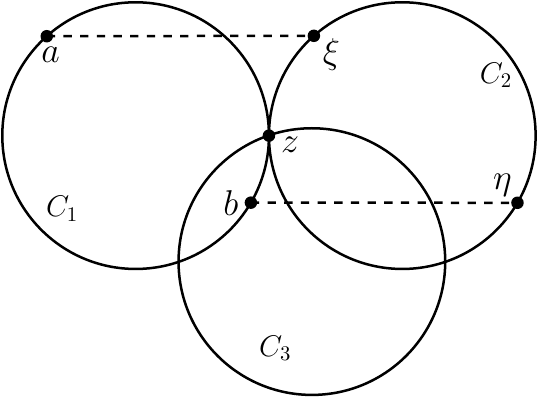}}
\hfill
\caption{(a) The numerator and denominator of the first fraction both vanish (the distance between the point $a$ and the 
center $c_2$ is $3$). (b) The exceptional situation in case N4D1 ($C_1$ is tangent to 
$C_2$ at $z$, and $\|a\xi\|=\|b\eta\|=2$).}
\label{N11D14}
\end{figure}

There are 16 cases of such a simultaneous vanishing of a numerator and a denominator.
The four numerators are 
(where $\tau_{u,C}$ denotes the tangent vector to the circle $C$ at the point $u$)
$$
(w-\xi)\cdot \tau_{\xi,C_2},\quad\quad
(z-w)\cdot \tau_{w,C_a},\quad\quad
(w'-z)\cdot \tau_{z,C_3},\quad\quad
(\eta-w')\cdot \tau_{w',C_b},
$$
and the four denominators are
$$
(w-\xi)\cdot \tau_{w,C_a},\quad\quad
(z-w)\cdot \tau_{z,C_3},\quad\quad
(w'-z)\cdot \tau_{w',C_b},\quad\quad
(\eta-w')\cdot \tau_{\eta,C_2}.
$$
As noted earlier, the vanishing of any of these eight expressions means that a corresponding
pair of points lie at distance $2$ from each other. Specifically, for the numerators, the
corresponding constraints are, respectively,
$$
\|wc_2\| = 2,\quad\quad
\|az\| = 2,\quad\quad
\|w'c_3\| = 2,\quad\quad
\|\eta b\| = 2,
$$
and for the denominators, the corresponding constraints are, respectively,
$$
\|a\xi\| = 2,\quad\quad
\|wc_3\| = 2,\quad\quad
\|zb\| = 2,\quad\quad
\|w'c_2\| = 2.
$$
Consider first the cases where the numerator and denominator of the same fraction both vanish.
Consider for specificity the first fraction, so we have
$$
(w-\xi)\cdot \tau_{\xi,C_2} =
(w-\xi)\cdot \tau_{w,C_a} = 0 ,\quad\quad\text{or}\quad\quad
\|a\xi\|=\|wc_2\|=2 . 
$$
In this case the four points $a$, $w$, $\xi$ and $c_2$ must be collinear, with
$\|a\xi\|=\|wc_2\|=2$. This is easily seen to imply that $\|ac_2\|=3$.
In this case, $a$ is an intersection point of $C_1$ with the circle of radius 
$3$ centered at $c_2$. Hence there are at most two choices for $a$, for a total of
$O(n)$ pairs $(a,b)$ that fall into this subcase. A similar argument applies when
the numerator and denominator of any of the three other fractions both vanish:
The corresponding constraints are $\|ac_3\|=3$, $\|bc_3\|=3$, and $\|bc_2\|=3$,
and in each of these cases it follows that there are at most two choices for either
$a$ or $b$, for a total of $O(n)$ pairs of these types. (Note that in these cases the 
curves $\gamma_{a,b}$ are singletons; see Figure~\ref{N11D14}(a).)

Therefore, in what follows,
we assume that the vanishing numerator and denominator belong to distinct fractions.
We use the mnemonic notation NtDs to mean that the numerator of the $t$-th fraction 
and the denominator of the $s$-th fraction both vanish, and consider the following
$12$ possible cases.

\noindent{\bf N1D2:}
Here we have $\|wc_3\|=\|wc_2\|=2$.
In this case, $w$ is an intersection point of the two circles of radius $2$ centered 
at $c_2$ and at $c_3$. Once $w$ is known, $a$ is an intersection point of $C_1$ with
the unit circle centered at $w$. Hence there are only $O(1)$ ways to choose $a$, for
a total of $O(n)$ ultra-degenerate pairs of this type.

\noindent{\bf N1D3:}
Here we have $\|wc_2\|=\|zb\|=2$.
We claim that each $b\in S_1$ can be coupled with only $O(1)$ $a$'s to form an
ultra-degenerate pair of this type. Indeed, given $b$, we find $z$, as an 
intersection point of $C_3$ and a circle of radius $2$ centered at $b$.
Given $z$, we find $w$, as an intersection point of a unit circle centered at
$z$ and a circle of radius $2$ centered at $c_2$. From $w$ we can find $a$, as 
an intersection point of $C_1$ and a unit circle centered at $w$. Altogether, there
are only $O(1)$ choices for $a$, as claimed.

\noindent{\bf N1D4:}
Here we have $\|wc_2\|=\|w'c_2\|=2$.
Here each $a\in S_1$ can be coupled with only $O(1)$ $b$'s. Indeed, given $a$,
we find $w$, as an intersection point of the unit circle centered at $a$ and 
a circle of radius $2$ centered at $c_2$. From $w$ we find $z$, as in the standard 
procedure, and from $z$ we find $w'$, as an intersection point of the unit circle 
centered at $z$ and a circle of radius $2$ centered at $c_2$. From $w'$ we find $b$, 
as an intersection point of $C_1$ with the unit circle centered at $w'$.
Altogether, there are only $O(1)$ choices for $b$, as claimed.

\noindent{\bf N2D1:}
Here we have $\|a\xi\|=\|az\|=2$.
Since $(a,\xi,z)$ span a unit circle, we must have $\xi=z$, which is thus an intersection
point of $C_2$ and $C_3$. This allows us to reconstruct $a$, essentially as in case N1D2.

\noindent{\bf N2D3:}
Here we have $\|az\|=\|zb\|=2$.
This case is easy: Given $a$, we find $z$, as an intersection point of $C_3$ and 
a circle of radius $2$ centered at $a$, and from $z$ we find $b$, as an intersection point 
of $C_1$ and a circle of radius $2$ centered at $z$. 

\noindent{\bf N2D4:}
Here we have $\|az\|=\|w'c_2\|=2$.
This case is symmetric to case N1D3, and is treated in a fully symmetric manner, starting
from $a$ and reconstructing $b$ in $O(1)$ ways.

\noindent{\bf N3D1:}
Here we have $\|w'c_3\|=\|a\xi\|=2$.
Given $a$, we find $\xi$, as an intersection point of $C_2$ and 
a circle of radius $2$ centered at $a$. From $a$ and $\xi$ we find $z$, as an intersection point 
of $C_3$ and the diametral (unit) circle determined by $a\xi$. From $z$ we find $w'$,
as the point that lies on the line through $c_3$ and $z$ at distance $1$ from $z$ (and $2$ from $c_3$).
From $w'$ we find $b$, as an intersection point of $C_1$ and the unit circle centered at $w'$.

\noindent{\bf N3D2:}
Here we have $\|w'c_3\|=\|wc_3\|=2$.
This case is treated exactly as case N1D4, except that here $c_3$ plays the role that was
played there by $c_2$.

\noindent{\bf N3D4:}
Here we have $\|w'c_3\|=\|w'c_2\|=2$, so this is a symmetric version of case N1D2,
with an essentially identical reconstruction process.

\noindent{\bf N4D1:}
Here we have $\|b\eta\|=\|a\xi\|=2$.
Given $a$, we find $\xi$, as an intersection point of $C_2$ with a circle of radius $2$ centered at $a$.
We then find $z$, as an intersection of $C_3$ with the diametral (unit) circle determined by $a\xi$.
In complete analogy with the treatment of case (iii) of the standard reconstruction process (at an 
extremal point of $\gamma'$), we note that the quadrilateral $R=c_1b\eta c_2$ has edges of fixed 
lengths, namely, $1$, $2$, $1$, and $\|c_1c_2\|$, and that it can flex around its fixed edge 
$c_1c_2$. As $R$ flexes, the midpoint of $b\eta$ traces an algebraic curve $\tau$ of some constant 
degree, and the intersection(s) of $\tau$ with the unit circle $C_z$ centered at $z$ gives us the center(s) 
$w'$ of the unit circle spanned by $(b,\eta,z)$, from which $b$ is readily obtained, as in several 
preceding cases. If $C_z$ and $\tau$ do not overlap, the number of intersection points between them 
is finite and bounded by a constant, and there are only $O(1)$ way to reconstruct $b$.

As in case (iii), the situation where $C_z$ and $\tau$ do overlap can happen only when $C_1$ and $C_2$
are tangent to each other, and $z$ is this tangency point (see Figure~\ref{N11D14}(b)). 
In this case it is possible to have a superlinear (in fact, any) number of pairs $(a,b)$, for which there exists 
$(\xi,\eta)\in C_2\times C_2$, such that $\|a\xi\|=\|b\eta\|=2$, and $(a,\xi,z)$, $(b,\eta,z)$ are unit triples. 
We therefore do not exclude those pairs as
ultra-degenerate. Instead, we claim that any such pair can be reconstructed from $\gamma'$, using
the reconstruction process described in Section~\ref{se:ov}. For this, we note that, for $(a,b)$ fixed, 
there exists at most one point $(\xi,\eta)$, 
such that $\|a\xi\|=\|b\eta\|=2$ (and $(a,\xi,z)$, $(b,\eta,z)$ are unit triples). Since, in the proof of 
Proposition~\ref{pr:ov}, we trace $\gamma'$ 
 from the point $q_0$, in both increasing and decreasing $v_x$-direction, we will not encounter 
this degeneracy in at least one of these traversals, and reach a locally $x$- or $y$-extremal point $(v_\xi,v_\eta)\in\gamma'$, 
from which $(a,b)$ can be reconstructed 
 (if not excluded in one of the other ultra-degenerate cases), in at most a constant number of ways.

\noindent{\bf N4D2:}
Here we have $\|b\eta\|=\|wc_3\|=2$.
This case is symmetric to case N3D1, and is treated in an analogous manner, starting from $b$
and reconstructing $a$ in $O(1)$ ways.

\noindent{\bf N4D3:}
Here we have $\|b\eta\|=\|bz\|=2$.
This is a symmetric variant of case N2D1, with an essentially identical treatment.

In summary, we have shown that the overall number of ultra-degenerate pairs is $O(n)$, as claimed.

\end{document}